\documentclass[a4]{elsart3-1}

\usepackage{amsmath}
\usepackage{amssymb}
\usepackage[english,francais]{babel}
\usepackage[latin1]{inputenc}
\usepackage[T1]{fontenc}
\usepackage[all]{xy}

\newtheorem{theorem}{Theorem}

\newtheorem{e-proposition}[theorem]{Proposition}
\newtheorem{corollary}[theorem]{Corollary}
\newtheorem{e-definition}[theorem]{Definition\rm}

\newtheorem{theoreme}{Th\'eor\`eme}

\newtheorem{proposition}[theoreme]{Proposition}
\newtheorem{corollaire}[theoreme]{Corollaire}

\renewcommand{\theequation}{\arabic{equation}}
\newenvironment{demonstration}
  {\setlength{\parindent}{0cm}\ \newline \textit{Démonstration }}
  {\setlength{\parindent}{5mm}}
\newcommand{\cqfd}{\hfill \ensuremath{\Box}}

\def\og{\leavevmode\raise.3ex\hbox{$\scriptscriptstyle\langle\!\langle$~}}
\def\fg{\leavevmode\raise.3ex\hbox{~$\!\scriptscriptstyle\,\rangle\!\rangle$}}

\makeatletter
\renewcommand\subsection{\@startsection{subsection}{2}{\z@}%
                                     {-3.25ex\@plus -1ex \@minus -.2ex}%
                                     {-3.00 mm}
                                     {\normalfont\bfseries}}
\makeatother
\renewcommand\subsection[1]{\vspace{.3cm}}

\newcommand{\C}{\mathbb{C}}
\newcommand{\Z}{\mathbb{Z}}
\renewcommand{\P}{\mathbb{P}}
\newcommand{\ig}{\mathbf{i}}
\newcommand{\jg}{\mathbf{j}}
\newcommand{\dg}{\mathbf{d}}
\DeclareSymbolFont{bbdold}{U}{bbold}{m}{n}
\DeclareSymbolFontAlphabet{\mathbbd}{bbdold}
\newcommand{\un}{\mathbbd 1}
\newcommand{\E}{\mathcal{E}^\dag}
\DeclareMathOperator{\add}{add}
\DeclareMathOperator{\Ext}{Ext}

\DeclareMathOperator{\md}{mod}
\DeclareMathOperator{\coker}{coker}
\makeatletter
\def\@fnsymbol#1{\ifcase#1\or \dagger\or \ddagger\else\@ctrerr\fi\relax} 
\makeatother
\long\def\symbolfootnote[#1]#2{\begingroup%
\def\thefootnote{\fnsymbol{footnote}}\footnote[#1]{#2}\endgroup}
\newcommand{\efface}[1]{}
\renewcommand{\epsilon}{\varepsilon}
\renewcommand{\phi}{\varphi}

\begin{document}

\begin{frontmatter}




%
\selectlanguage{francais}
\title{Algèbres amassées et algèbres préprojectives : le cas non simplement lacé}

\vspace{-2.6cm}
\selectlanguage{english}
\title{Cluster algebras and preprojective algebras : the non simply-laced case}



\author[LD]{Laurent Demonet}
\ead{Laurent.Demonet@normalesup.org}

\address[LD]{LMNO, Université de Caen, Esplanade de la Paix, 14000 Caen}

\begin{abstract}
We generalize to the non simply-laced case results of Gei\ss, Leclerc and Schröer about the cluster structure of the coordinate ring of the maximal unipotent subgroups of simple Lie groups. In this way, cluster structures in the non simply-laced case can be seen as projections of cluster structures in the simply-laced case. This allows us to prove that cluster monomials are linearly independent in the non simply-laced case.
\efface{\it To cite this article: L. Demonet, C. R.
Acad. Sci. Paris, Ser. I 336 (2003).}

\vskip 0.5\baselineskip

\selectlanguage{francais}
\noindent{\bf R\'esum\'e}
\vskip 0.5\baselineskip
\noindent
On généralise au cas non simplement lacé des résultats de Gei\ss, Leclerc et Schröer sur les structures amassées des algèbres de fonctions sur les sous-groupes unipotents maximaux des groupes de Lie simples. Cela permet en particulier de voir les structures amassées dans le cas non simplement lacé comme projections des structures amassées dans le cas simplement lacé. Cela permet aussi de montrer la liberté des monômes d'amas dans le cas non simplement lacé. {\it Pour citer cet article~: L. Demonet, C. R.
Acad. Sci. Paris, Ser. I 336 (2003).}

\end{abstract}
\end{frontmatter}

\selectlanguage{english}
\section*{Abridged English version}

Fomin and Zelevinsky \cite{FZ1} introduced cluster algebras and constucted with Berenstein \cite{BFZ} a cluster structure on the function algebras $\C[N]$ where $N$ is a maximal unipotent subgroup of a simple Lie group $G$. This structure was interpreted in the simply-laced case ($A$, $D$, $E$) by Gei\ss, Leclerc and Schröer in terms of rigid representations of preprojective algebras \cite{GLS3}. Our aim is to extend these results to the non simply-laced case ($B$, $C$, $F$ and $G$). 

Let $D$ be a simply-laced Dynkin diagram and $a$ be an admissible automorphism of $D$. There are two natural ways to construct from $D$ and $a$ another Dynkin diagram : $D'$ defined by Lusztig \cite{L3} and $D''$ defined by Kac \cite{Kac}. All the non simply-laced Dynkin diagrams can be obtained in each of these two ways. For example, if $D=A_{2n-1}$ and $a$ is the order $2$ diagram automorphism, we have $D' = B_n$ and $D'' = C_n$.

Let $G$ be the connected and simply connected Lie group of Dynkin diagram $D$ and $G''$ be the connected and simply connected Lie group of Dynkin diagram $D''$. Let $N$ and $N''$ be their maximal unipotent subgroups. Let  $\mathfrak{n}$ and $\mathfrak{n}''$ be the Lie algebras of $N$ and $N''$. There is a natural injective morphism from $\mathfrak{n}''$ into $\mathfrak{n}$ \cite[7.9]{Kac}. Therefore $U(\mathfrak{n}'')$ can be identified to a subalgebra of $U(\mathfrak{n})$. As $U(\mathfrak{n})^* \simeq \C[N]$ and $U(\mathfrak{n''})^* \simeq \C[N'']$ as Hopf algebras, this gives a surjection $p:\C[N] \twoheadrightarrow \C[N'']$.

Let $A$ be a ring and $\Gamma$ be a group acting on $A$. The additive group $A \otimes_\Z \Z[\Gamma]$ can be endowed with a ring structure by defining $(a_1 \otimes g_1) (a_2 \otimes g_2) = (a_1 g_1(a_2) \otimes g_1 g_2)$ for all $(a_1, g_1), (a_2, g_2) \in A \times \Gamma$. This ring, called a skew-group ring of $\Gamma$ in \cite{RR}, will be denoted by $A \Gamma$. If $A$ is a $k$-algebra, it induces a $k$-algebra structure on $A\Gamma$.

Let $\Delta$ be a Dynkin diagram with set of vertices $I$. Let $\Delta_1$, $\Delta_2$ be simply-laced Dynkin diagrams with sets of vertices $I_1$, $I_2$ and $a_1$, $a_2$ be admissible automorphisms of $\Delta_1$ and $\Delta_2$ such that $\Delta_1' = \Delta_2'' = \Delta$. Let $n$ be the common order of $a_1$ and $a_2$. Let $\Lambda_1$ and $\Lambda_2$ be the preprojective algebras of some orientations of $\Delta_1$ and $\Delta_2$ compatible with $a_1$ and $a_2$ (see \cite{GLS1}). Let $\Gamma_1$ be the automorphism group of $\Delta_1$ generated by $a_1$. Using methods of \cite{RR}, one proves
\begin{e-proposition}
 There is an equivalence of categories $\Phi : \md \Lambda_1\Gamma_1 \rightarrow \md \Lambda_2$.
\end{e-proposition}

Let $F:\md \Lambda_1 \Gamma_1 \rightarrow \md \Lambda_1$ be the forgetful functor. For  $M \in \md \Lambda_2$, set $S(M) = \bigoplus_{k = 0}^{n-1} a_2^k(M)$ and $\Sigma(M)$ the largest basic submodule of $S(M)$.

Two $\Lambda_2$-modules $M$ and $N$ will be called $a_2$-isomorphic if $S(M) \simeq S(N)$.

A $\Lambda_1 \Gamma_1$-module $M$ will be called $F$-rigid (resp. $F$-basic) if $F(M)$ is rigid (resp. basic). A $\Lambda_2$-module $M$ will be called $a_2$-rigid (resp. $a_2$-basic) if $S(M)$ is rigid (resp. if it has not two direct summands which are $a_2$-isomorphic). An $F$-rigid $\Lambda_1 \Gamma_1$-module (resp. $a_2$-rigid $\Lambda_2$-module) $M$ will be called maximal $F$-rigid (resp. $a_2$-rigid) if for every $\Lambda_1 \Gamma_1$-module (resp. $\Lambda_2$-module) $X$, if $M \oplus X$ is $F$-rigid (resp. $a_2$-rigid) then $F(X)$ (resp. $S(X)$) is a direct summand of $F(M)$ (resp. $S(M)$).

For $\ig \in I = I_2/a_2$, define $Q_{\ig} = \bigoplus_{i \in \ig} Q_i$ where $Q_i$ is the injective envelope of the simple $\Lambda_2$-module $S_i$ and $\E_\ig = \prod_{i \in \ig} \E_i$ where the $\E_i$ are defined in \cite{GLS6}. Let $\ig_1 \ig_2 \dots \ig_k$ be a reduced expression of the longest element of the Weyl group associated to $\Delta$. One can define the maximal rigid basic $\Lambda_2$-module
 $$T_{\ig_1 \ig_2 \dots \ig_k} = \bigoplus_{\ell = 1}^k \E_{\ig_1} \E_{\ig_2} \dots \E_{\ig_\ell} \left( Q_{\ig_\ell}\right) \oplus \bigoplus_{\ig \in I} Q_\ig.$$

\begin{e-proposition}
 \begin{enumerate}
  \item Let $X, X' \in \md \Lambda_1 \Gamma_1$. Then $F(X) \simeq F(X')$ if and only if $S(\Phi(X)) \simeq S(\Phi(X'))$.
  \item A $\Lambda_1 \Gamma_1$-module $X$ is $F$-rigid (resp. $F$-basic) if and only if $\Phi(X)$ is $a_2$-rigid (resp. $a_2$-basic).
  \item Let $T_2 \in \md \Lambda_2$ such that $\Sigma(T_2) = T_{\ig_1 \ig_2 \dots \ig_k}$. It is maximal $a_2$-rigid and has $r$ non $a_2$-isomorphic indecomposable summands where $r$ is the number of positive roots of $\Delta$. 
  \item The $a_2$-rigid $\Lambda_2$-modules have at most $r$ non $a_2$-isomorphic summands. The maximal ones have exactly $r$ non $a_2$-isomorphic summands.
  \item A $\Lambda_2$-module is maximal $a_2$-rigid if and only if $\Sigma(X)$ is maximal rigid. Moreover, each maximal basic rigid $a_2$-stable $\Lambda_2$-module has at least one preimage by $\Sigma$.
  \end{enumerate}
\end{e-proposition}

Let $T = T_0 \oplus X$ be a maximal $a_2$-rigid $a_2$-basic $\Lambda_2$-module, such that $X$ is a non projective indecomposable summand. If $f$ is a minimal left $\add(S(T_0))$-approximation of $X$ then $f$ is injective and we will denote $\nu_X(T) = T_0 \oplus Y$ where $Y = \coker f$.
\begin{e-proposition}
 \begin{enumerate}
  \item The module $Y$ is indecomposable and $\nu_X(T)$ is maximal $a_2$-rigid $a_2$-basic.
  \item $\nu_X(\nu_Y(T))=T$.
  \item $\Sigma(\nu_X(T))$ only depends on $\Sigma(T)$ and $\Sigma(X)$.
 \end{enumerate}
\end{e-proposition}

\efface{
If $\mu$ denotes the mutation operator on the maximal rigid basic $\Lambda_2$-modules \cite{GLS3} and $\tilde X$ denotes the set of indecomposable summands of $\Sigma(X)$, we have the following result :
\begin{e-proposition}
 $$\Sigma \left(\nu_X(T)\right) = \prod_{X' \in \tilde X} \mu_{X'} \left(\Sigma(T)\right)$$
\end{e-proposition}
}
Let $G_2$ and $G$ be the connected and simply connected Lie groups of type $\Delta_2$ and $\Delta$. Let $N_2$ and $N$ be their maximal unipotent subgroups. Each $\Lambda_2$-module $X$ gives rise to $\phi_{X} \in \C[N_2]$ (\cite{L1} and \cite{GLS5}). Let $\psi_X = p(\phi_{X})$, where $p:\C[N_2] \rightarrow \C[N]$ is defined as above.

\begin{e-proposition}
 \begin{enumerate}
  \item If $X \in \md \Lambda_2$ then $\psi_X$ only depends on $\Sigma(X)$. 
  \item Let $T$ be a maximal $a_2$-rigid $\Lambda_2$-module such that $\Sigma(T) = T_{\ig_1 \ig_2 \dots \ig_k}$. Then the $\psi_X$'s where $X$ runs over the indecomposable summands of $T$ form an initial seed of \cite{BFZ} for the cluster algebra $\C[N]$. 
  \item Let $\mathcal E$ be the set of the maximal $a_2$-rigid $a_2$-basic $\Lambda_2$-modules that can be reached using mutations of type $\nu$ from $T$. There is a bijection between the isomorphism classes of $\Sigma(\mathcal E)$ and the clusters of $\C[N]$.
 \end{enumerate}
\end{e-proposition}

 All previous results in terms of $\Lambda_2$ and $\Sigma$ can be translated in terms of $\Lambda_1 \Gamma_1$ and $F$.

\begin{corollary}
 The cluster variables of $\C[N]$ are $f \in \C[N]$ such that $p^{-1}(f)$ is a collection of cluster variables of the same cluster of $\C[N_2]$. The preimage by $p$ of the clusters of $\C[N]$ are clusters of $\C[N_2]$.
\end{corollary}

This allows to prove a conjecture of Fomin and Zelevinsky which was proved in the simply-laced case by Gei\ss, Leclerc et Schröer :
\begin{theorem}
 The cluster monomials of $\C[N]$ are linearly independent.
\end{theorem}%
%
%
%
%
%
%
\selectlanguage{francais}%
%
\section{Définitions et notations}
\subsection{}

Soit un diagramme de Dynkin $D$ simplement lacé d'ensemble de sommets $I$ et $a$ un automorphisme admissible de $D$ (c'est-à-dire que deux sommets d'une même orbite ne sont jamais reliés par une flèche). Soit $C$ la matrice de Cartan de $D$. Il y a deux façons naturelles de construire une nouvelle matrice de Cartan indexée par l'ensemble $I/a$ des orbites :
$$C'_{\ig \jg} = \frac{1}{\# \ig} \sum_{(i,j) \in \ig \times \jg} C_{ij} \qquad\text{ et }\qquad C''_{\ig \jg} = \frac{1}{\# \jg} \sum_{(i,j) \in \ig \times \jg} C_{ij}$$
où $\ig, \jg \in I/a$.
On notera $D'$ et $D''$ les diagrammes de Dynkin correspondants. Par exemple, si $D = A_{2n-1}$ et $a$ est l'unique automorphisme d'ordre $2$, alors $D' = B_n$ et $D''=C_n$. On peut construire tous les diagrammes de Dynkin non simplement lacés de chacune de ces deux façons.
\efface{\begin{center}
 \begin{tabular}{|l|c|c|c|c|}
  \hline 
   type de $C$ & $A_{2n-1}$ & $D_{n+1} $ & $D_4$ & $E_6$ \\
  \hline
   ordre de $a$ & $2$ & $2$ & $3$ & $2$ \\
  \hline
   type de $C'$ & $B_n$ & $C_n$ & $G_2$ & $F_4$ \\
  \hline
   type de $C''$ & $C_n$ & $B_n$ & $G_2$ & $F_4$ \\
  \hline
 \end{tabular}
\end{center}}

\subsection{}
\label{lie}
Soit $G$ le groupe de Lie connexe et simplement connexe de diagramme de Dynkin $D$ et $G''$ le groupe de Lie connexe et simplement connexe de diagramme de Dynkin $D''$. Soient $N$ et $N''$ des sous-groupes unipotents maximaux de $G$ et de $G''$. Soient aussi $\mathfrak{n}$ et $\mathfrak{n}''$ leurs algèbres de Lie. Il y a un unique morphime injectif $\iota$ de $\mathfrak{n}''$ dans $\mathfrak{n}$ vérifiant $\iota(e''_{\ig}) = \sum_{i \in \ig} e_i$
où les $e''_{\ig}$ désignent les générateurs de Chevalley de $\mathfrak{n}''$ et les $e_i$ ceux de $\mathfrak{n}$ (voir \cite[7.9]{Kac}). On peut donc identifier $U(\mathfrak{n}'')$ à une sous-algèbre de $U(\mathfrak{n})$.
\efface{ Du couplage 
\begin{align*}
 U(\mathfrak{n}) \times \C[N] &\rightarrow \C \\
 (e_{i_1} e_{i_2} \dots e_{i_k}, f) &\mapsto \left.\frac{\partial}{\partial t_1} \frac{\partial}{\partial t_2} \dots \frac{\partial}{\partial t_k} f\left(\exp(t_1 e_{i_1})\exp(t_2 e_{i_2})\dots \exp(t_k e_{i_k})\right)\right|_{t_1 = t_2 = \dots = t_k = 0}
\end{align*}
on déduit un isomorphisme $U(\mathfrak{n})^* \simeq \C[N]$.} Rappelons que $U(\mathfrak{n})^* \simeq \C[N]$ et $U(\mathfrak{n''})^* \simeq \C[N'']$ comme algèbres de Hopf. Finalement, on a une surjection $p:\C[N] \twoheadrightarrow \C[N'']$ qui correspond au fait que $N''$ est un sous-groupe algébrique de $N$.

\subsection{}
Si $A$ est un anneau et $\Gamma$ est un groupe agissant sur $A$, on notera $A \Gamma$ l'anneau dont le groupe additif est $A \otimes_\Z \Z[\Gamma]$ et la multiplication prolonge $(a_1 \otimes g_1) (a_2 \otimes g_2) = (a_1 g_1(a_2) \otimes g_1 g_2)$ pour $(a_1, g_1), (a_2, g_2) \in A \times \Gamma$. Si $A$ a une structure de $k$-algèbre, elle induit une structure de $k$-algèbre sur $A\Gamma$. Pour plus de détails sur cette construction, voir \cite{RR}.

\section{Modules rigides symétriques et mutations}
\subsection{}
Soit $\Delta$ un diagramme de Dynkin quelconque. Soient $\Delta_1$ (resp. $\Delta_2$) le diagramme de Dynkin simplement lacé et $a_1$ (resp. $a_2$) l'automorphisme de $\Delta_1$ (resp. $\Delta_2$) vérifiant $\Delta'_1 = \Delta$ (resp. $\Delta''_2 = \Delta$). Notons $n$ l'ordre commun de $a_1$ et de $a_2$. Soient aussi $Q_1 = (I_1, F_1)$ et $Q_2 = (I_2, F_2)$ des orientations de $\Delta_1$ et $\Delta_2$ compatibles respectivement avec $a_1$ et $a_2$. Notons enfin $\Lambda_1$ et $\Lambda_2$ les algèbres préprojectives de $Q_1$ et $Q_2$ (voir \cite{GLS1} ou \cite{L1}). \efface{Si $\ig$ est un sommet de $\Delta$, l'orbite correspondante dans $\Delta_1$ sera notée $\ig_1$ et l'orbite correspondante dans $\Delta_2$ sera notée $\ig_2$. Il est alors aisé de remarquer que $\# \ig_1 \# \ig_2 = n$ pour tout sommet $\ig$ de $\Delta$.}

Soit $\Gamma_1$ le groupe d'automorphismes de $\Delta_1$ engendré par $a_1$. Ce groupe agit sur $\Lambda_1$ donc on peut former $\Lambda_1 \Gamma_1$.
\begin{proposition}
 Il existe une équivalence de catégories $\Phi : \md \Lambda_1\Gamma_1 \rightarrow \md \Lambda_2$.
\end{proposition}

On démontre ce résultat à l'aide des méthodes présentées en \cite{RR}.

\subsection{}
Notons $F:\md \Lambda_1 \Gamma_1 \rightarrow \md \Lambda_1$ le foncteur d'oubli. Pour tout $M \in \md \Lambda_2$, notons $S(M) = \bigoplus_{k = 0}^{n-1} a_2^k(M)$ et $\Sigma(M)$ le plus grand sous-module basique de $S(M)$.

Deux $\Lambda_2$-modules $M$ et $N$ seront dit $a_2$-isomorphes si $S(M) \simeq S(N)$.

Un $\Lambda_1 \Gamma_1$-module (resp. $\Lambda_2$-module) $M$ sera appelé $F$-rigide (resp. $a_2$-rigide) si $F(M)$ (resp. $S(M)$) est rigide. Il sera dit $F$-rigide maximal (resp. $a_2$-rigide maximal) si pour tout $\Lambda_1 \Gamma_1$-module (resp. $\Lambda_2$-module) indécomposable $X$, si $M \oplus X$ est $F$-rigide (resp $a_2$-rigide) alors $F(X)$ (resp. $S(X)$) est un facteur indécomposable de $F(M)$ (resp. $S(M)$).

Un $\Lambda_1 \Gamma_1$-module $M$ sera appelé $F$-basique si $F(M)$ est basique. Un $\Lambda_2$-module $M$ sera appelé $a_2$-basique si il n'a pas deux facteurs directs $a_2$-isomorphes.

Pour $\ig \in I = I_2/a_2$, on définit $Q_{\ig} = \bigoplus_{i \in \ig} Q_i$ où les $Q_i$ sont les enveloppes injectives des $\Lambda_2$-modules simples $S_i$ et $\E_\ig = \prod_{i \in \ig} \E_i$ où les $\E_i$ sont les foncteurs définis dans \cite{GLS6} ; $\E_\ig$ est bien défini car les $\E_i$ intervenant commutent. Soit maintenant $\ig_1 \ig_2 \dots \ig_k$ une expression réduite du mot de plus grande longueur du groupe de Weyl correspondant à $\Delta$ ; on peut lui associer le $\Lambda_2$-module
 $$T_{\ig_1 \ig_2 \dots \ig_k} = \bigoplus_{\ell = 1}^k \E_{\ig_1} \E_{\ig_2} \dots \E_{\ig_\ell} \left( Q_{\ig_\ell}\right) \oplus \bigoplus_{\ig \in I} Q_\ig.$$
 C'est un module rigide basique maximal.
\begin{proposition}
 \begin{enumerate}
  \item Soient $X, X' \in \md \Lambda_1 \Gamma_1$. Alors $F(X) \simeq F(X')$ si et seulement si $S(\Phi(X)) \simeq S(\Phi(X'))$.
  \item Un $\Lambda_1 \Gamma_1$-module $X$ est $F$-rigide (resp. $F$-basique) si et seulement si $\Phi(X)$ est $a_2$-rigide (resp. $a_2$-basique).
   \item Soit $T_2 \in \md \Lambda_2$ $a_2$-basique tel que $\Sigma(T_2) = T_{\ig_1 \ig_2 \dots \ig_k}$. Ce module est $a_2$-rigide maximal et il a $r$ facteurs directs indécomposables non $a_2$-isomorphes où $r$ est le nombre de racines postives du système associé à $\Delta$. 
  \item Les $\Lambda_2$-modules $a_2$-rigides ont au plus $r$ facteurs directs non $a_2$-isomorphes. De plus, ceux qui sont maximaux ont exactement $r$ facteurs directs non $a_2$-isomorphes. 
  \item Un $\Lambda_2$-module $X$ est $a_2$-rigide maximal si et seulement si $\Sigma(X)$ est rigide maximal. De plus, tout $\Lambda_2$-module rigide basique maximal $a_2$-stable admet au moins un antécédent par $\Sigma$.
  \efface{\item Un module $\Lambda_2$-$X$ est $a_2$-rigide $a_2$-basique maximal si et seulement si $\Sigma(X)$ est rigide basique maximal. De plus, tout $\Lambda_2$-module rigide basique maximal $a_2$-stable admet au moins un antécédent par $\Sigma$.}
 \end{enumerate}
\end{proposition}

\subsection{}
Soit $T = T_0 \oplus X$ un $\Lambda_2$-module $a_2$-rigide maximal et $a_2$-basique, tel que $X$ soit un facteur indécomposable non projectif. Si $f$ est une $\add(S(T_0))$-approximation minimale à gauche de $X$, alors $f$ est injective et l'on notera $\nu_X(T) = T_0 \oplus Y$ où $Y$ est le conoyau de $f$.
\begin{proposition}
 \begin{enumerate}
  \item Le module $Y$ est indécomposable et $\nu_X(T)$ est $a_2$-rigide maximal et $a_2$-basique.
  \item $\nu_X(\nu_Y(T))=T$.
  \item $\Sigma(\nu_X(T))$ ne dépend que de $\Sigma(T)$ et $\Sigma(X)$.
  \efface{\item On peut transporter à travers $\Phi$ l'opération de mutation $\nu$ en une opération de mutation sur les modules $F$-rigides et $F$-basiques maximaux de $\Lambda_1\Gamma_1$. On la notera aussi $\nu$.
  \item Si $T'$ est un $\Lambda_1 \Gamma_1$-module $F$-rigide et $F$-basique maximal et $X'$ un de ses facteurs non projectifs alors $F(\nu_{X'}(T'))$ ne dépend que de $F(T')$ et $F(X')$.}
 \end{enumerate}
\end{proposition}

Le module $\Sigma(T)$ est un $\Lambda_2$-module rigide basique maximal invariant par l'action de $a_2$. On peut considérer, suivant \cite{GLS3}, pour tout $\Lambda_2$-module rigide basique maximal $T'$ et pour tout facteur indécomposable $X'$ de $T'$ le module $\mu_{X'}(T')$ obtenu par mutation. Notons $\tilde X$ l'ensemble des facteurs indécomposables de $\Sigma(X)$. On montre que les mutations $\nu$ s'expriment de la manière suivante à partir des mutations $\mu$ :
 $$\Sigma \left(\nu_X(T)\right) = \prod_{X' \in \tilde X} \mu_{X'} \left(\Sigma(T)\right)$$
 où le second membre ne dépend pas de l'ordre des mutations.
 
\section{Algèbres amassées}

\subsection{}
Soient $G_i$ et $G$ les groupes de Lie connexes et simplement connexes de type $\Delta_i$ et $\Delta$ ($i \in \{1,2\}$). Soient $N_i$ et $N$ leurs sous groupes unipotents maximaux et $\mathfrak{n}_i$ et $\mathfrak{n}$ les algèbres de Lie de $N_i$ et $N$. \`A chaque $\Lambda_2$-module $X$, on peut associer une fonction $\delta_{X}$ de $U(\mathfrak{n}_2)^*$ (voir \cite{L1}) et donc une fonction $\phi_{X}$ de $\C[N_2]$ (voir \ref{lie}). Notons $\psi_X = p(\phi_{X})$, où $p$ est définie en \ref{lie}. Si $X'$ est un $\Lambda_1 \Gamma_1$-module, on notera aussi $\psi_{X'} = \psi_{\Phi(X')}$.

\begin{proposition}
 \begin{enumerate}
  \item Pour tout $X \in \md \Lambda_2$, $\psi_X = \psi_{a_2(X)}$. Autrement dit $\psi_X$ ne dépend que de $\Sigma(X)$. 
  \item Soit $T$ un module $a_2$-rigide maximal tel que $\Sigma(T) = T_{\ig_1 \ig_2 \dots \ig_k}$. Alors l'ensemble des $\psi_X$ où $X$ est un facteur indécomposable de $T$ forme l'une des graines initiales de \cite{BFZ} pour l'algèbre amassée $\C[N]$. Cette graine initiale ne dépend que du mot $\ig_1 \ig_2 \dots \ig_k$ et on obtient toutes les graines initiales de \cite{BFZ} de cette façon.
  \item Soit $\mathcal{E}$ l'ensemble des $\Lambda_2$-modules $a_2$-rigides maximaux que l'on peut atteindre par la mutation $\nu$ à partir de $T$. Il y a une bijection entre les classes d'isomorphisme de $\Sigma(\mathcal E)$ et les amas de l'algèbre amassée $\C[N]$.
  \efface{\item Soit $T'$ un $\Lambda_1 \Gamma_1$-module tel que $\Phi(T') = T$. Soit $E_1$ l'ensemble des $\Lambda_1 \Gamma_1$-modules $F$-rigides $F$-basiques maximaux que l'on peut atteindre par la mutation $\nu$ à partir de $T'$. Il y a une bijection entre les classes d'isomorphisme de $F(E_1)$ et les amas de l'algèbre amassée $\C[N]$.}
 \end{enumerate}
\end{proposition}

Tous les résultats précédents concernant $\Lambda_2$ et $\Sigma$ se traduisent par des résultats équivalents concernant $\Lambda_1 \Gamma_1$ et $F$.

\efface{
\subsection{}
Plus précisément, on a la formule de multiplication suivante, qui est traduction immédiate de celle de Gei\ss, Leclerc et Schröer (voir \cite{GLS5}) :
\begin{proposition}
 Si $X$ et $Y$ sont deux $\Lambda_2$-modules rigides symétriques,
 $$\chi(\P \Ext^1(X, Y)) \delta_X \delta_Y = \sum_{Z \in E_{X,Y}} \chi(\P \Ext^1_Z(X, Y)) \delta_Z + \sum_{Z \in E_{Y,X}} \chi(\P \Ext^1_Z(Y, X)) \delta_Z$$
 où $E_{X, Y}$ désigne l'ensemble des classes d'extensions non scindées de $X$ par $Y$ et $\Ext^1_Z(X, Y)$ désigne la sous variété de $\Ext^1(X, Y)$ des extensions isomorphes à $Z$. Par ailleurs, $\chi$ désigne la caractéristique d'Euler-Poincaré.
\end{proposition}

Cette formule doit en fait être interprétée comme un ensemble de formules de multiplication puisque $\delta_{a_2(X)} = \delta_X$ alors que la formule obtenue en remplaçant $X$ par $a_2(X)$ n'est en général pas la même. Dans le cas où $Y$ est obtenue à partir de $X$ par mutation, on ne trouve qu'une formule non triviale (c'est à dire où $\Ext^1(X', Y') \neq 0$) en faisant varier $X'$ et $Y'$ tels que $X'$ soit dans la même orbite que $X$ par $a_2$ et $Y'$ dans la même orbite que $Y$. Cette formule est la formule de mutation des variables d'amas.
}

\subsection{}
On obtient le corollaire suivant :
\begin{corollaire}
 Les variables d'amas de $\C[N]$ sont des fonctions $f$ sur $N$ telles que $p^{-1}(f)$ soit consituée de variables d'amas appartenant toutes à un même amas de $\C[N_2]$. L'image réciproque par $p$ d'un amas de $\C[N]$ est un amas de $\C[N_2]$.
 
 De plus, si $\mathbf x$ est un amas de $\C[N]$, et $x \in \mathbf x$ n'est pas dans l'anneau des coefficients de l'algèbre amassée (autrement dit s'il ne provient pas d'un $\Lambda_2$-module projectif), 
 $$p^{-1}(\nu_x (\mathbf x)) = \prod_{p(x_0) = x} \mu_{x_0} \left(p^{-1}(\mathbf x)\right)$$
 la composition étant bien définie car les $\mu_{x_0}$ commutent dans ce cas.
\end{corollaire}

Ce corollaire peut aussi se démontrer directement avec les mêmes techniques que dans \cite{GD}.

\efface{
Une autre manière d'exprimer la dernière partie de ce corollaire est de dire que si $B^2$ est la matrice d'échange de l'amas $p^{-1}(\mathbf x)$ (voir \cite{FZ1}), alors la matrice d'échange $B$ de l'amas $\mathbf x$ est celle vérifiant, pour tout $x \in \mathbf x$ qui n'est pas dans l'anneau des coefficients et $y \in \mathbf x$,
$$B_{yx} = \sum_{p(y_0) = y} B^2_{y_0 x_0}$$
ceci ne dépendant pas du choix de $x_0 \in p^{-1}(x)$.

De la même façon, on peut montrer, si $B^1$ est la matrice d'échange de l'amas correspondant à l'amas $p'^{-1}(\mathbf x)$ de $\C[N^1]$ (à travers le foncteur entre $\Lambda_2$ et $\Lambda_1 \Gamma_1$ et puis par projection sur $\Lambda_1$) que pour tout $x \in \mathbf x$ qui n'est pas dans l'anneau des coefficients et tout $y \in \mathbf x$,
$$B_{yx} = \sum_{p'(x_0) = x} B^1_{y_0 x_0}$$
ceci ne dépendant pas du choix de $y_0 \in p'^{-1}(y)$ et $p'$ étant l'application qui à un module rigide sur $\Lambda_1$ associe une classe modulo $a_2$ de $\Lambda_2$-modules rigides symétriques.
}

\subsection{}
Si $M \in \md \Lambda_1$ est stable par $a_1$ (i.e. vérifie $a_1 M \simeq M$) et $\ig \in I = I_1/a_1$, on note $k_{\ig}(M)$ la dimension du sous-module maximal de $M$ supporté par $\ig$. On note $\mathcal{R}$ l'ensemble des classes d'isomorphisme de $\Lambda_1$-modules rigides stables par $a_1$, et si $\dg$ est un vecteur dimension stable par $a_1$, $\mathcal{R}_\dg$ désigne l'ensemble des éléments de $\mathcal{R}$ de dimension $\dg$ et $\mathcal{R}_{\dg, \ig, k} = \{M \in \mathcal{R}_\dg \,|\, k_{\ig}(M) = k\}$. On montre que $\E_\ig$ induit une injection de $\mathcal{R}_{\dg, \ig, k}$ dans $\mathcal{R}_{\dg - k \ig, \ig, 0}$.

Pour tout $M \in \md \Lambda_1$ stable par $a_1$, notons $\epsilon_M$ l'élément de $U(\mathfrak n)^*$ correspondant à $\psi_M$.

\begin{theoreme}
 \begin{enumerate}
  \item Pour tout $M \in \mathcal{R}_\dg$, il existe $f_M \in U(\mathfrak n)$ homogène de degré $\dg$ tel que pour tout $N \in \mathcal{R}$, $\epsilon_N (f_M) = \un_{M=N}$.
  \item Les monômes d'amas (c'est-à-dire les produits de variables d'amas d'un même amas) de $\C[N]$ sont linéairement indépendants.
 \end{enumerate}
\end{theoreme}

\begin{demonstration}
 \begin{enumerate}
  \item On construit $f_M$ par récurrence sur $\dg$. Pour $\dg = 0$, $f_0 = 1$. Supposons $f_{M'}$ construit pour tout $M'$ de dimension strictement inférieure à $\dg$. Comme $M$ est nilpotente, il existe $\ig$ tel que $k_{\ig}(M) > 0$. Raisonnons donc par récurrence descendante sur $k_{\ig}(M)$. Supposons le résultat montré pour tous les $M'$ de dimension $\dg$ et tels que $k_{\ig}(M) < k_{\ig}(M') \leq \dg_\ig$. Soit $f_0 = f_{\E_\ig(M)} e_\ig^{k_{\ig}(M)}$. Alors
  $$f_M = f_0 - \sum_{N \in \mathcal{R}_\dg \,|\, k_{\ig}(N) > k_{\ig}(M)} \epsilon_N (f_0) f_N.$$
  convient.
  \item On déduit du point précédent que les $\epsilon_M$ sont linéairement indépendants, donc les $\psi_M$, et on conclut en observant que les monômes d'amas sont exactement les $\psi_M$, où $M$ parcourt $\mathcal{R}$. \cqfd
 \end{enumerate}

\end{demonstration}

Le second point est une conjecture de Fomin et Zelevinsky, qui avait été montré dans le cas simplement lacé par Gei\ss, Leclerc et Schröer.
\efface{

\section*{Exemple}

Ici, plaçons nous dans le cas où $D_0 = A_5 = \xymatrix{3 \ar@{-}[r] & 2 \ar@{-}[r] & 1 \ar@{-}[r] & 2' \ar@{-}[r] & 3'}$ et $D = C_3 = \xymatrix{1 \ar@{=>}[r] & 2 \ar@{-}[r] & 3}$.

}

\efface{
\section*{Remerciements}
Je remercie B. Leclerc pour les discussions intéressantes et fructueuses à ce sujet et pour la relecture attentive de cette note. Je remercie aussi C. Gei\ss\ et I. Reiten pour des discussions intéressantes à propos de ces constructions.
}


\begin{thebibliography}{00}

%
  \bibitem{BFZ} A. Berenstein, S. Fomin, A. Zelevinsky, \textit{Cluster algebras III: Upper bounds and double Bruhat cells}, Duke Math. J.  126  (2005),  no. 1, 1--52. 
  \bibitem{GD} G. Dupont, \textit{An approach to non simply laced cluster algebras},  arXiv: math/0512043.
  \bibitem{FZ1} S. Fomin, A. Zelevinsky, \textit{Cluster algebras I: Foundations}, J. Amer. Math. Soc.  15  (2002),  no. 2, 497--529
  \bibitem{GLS1} C. Gei\ss, B. Leclerc, J. Schr{\"o}er,
  \textit{Semicanonical bases and preprojective algebras},  Ann. Sci. École Norm. Sup. (4)  38  (2005),  no. 2, 193--253.
  \bibitem{GLS2} C. Gei\ss, B. Leclerc, J. Schr{\"o}er,
  \textit{Verma modules and preprojective algebras},  Nagoya Math. J.  182  (2006), 241--258.
  \bibitem{GLS3} C. Gei\ss, B. Leclerc, J. Schr{\"o}er,
  \textit{Rigid modules over preprojective algebras},  Invent. Math.  165  (2006),  no. 3, 589--632. 
  \bibitem{GLS4} C. Gei\ss, B. Leclerc, J. Schr{\"o}er,
  \textit{Auslander algebras and initial seeds for cluster algebras},  J. London Math. Soc. 75  (2007),  no. 3, 718--740.
  \bibitem{GLS5} C. Gei\ss, B. Leclerc, J. Schr{\"o}er,
  \textit{Semicanonical bases and preprojective algebras II: A multiplication formula},  Compositio Math. 143  (2007),  no. 5, 1313--1334.
  \bibitem{GLS6} C. Gei\ss, B. Leclerc, J. Schr{\"o}er,
  \textit{Partial flag varieties and preprojective algebras}, arXiv: math/0609138. Ann. Inst. Fourier, à paraître.
  \bibitem{GLS7} C. Gei\ss, B. Leclerc, J. Schr{\"o}er,
  \textit{Rigid modules over preprojective algebras II: The Kac-Moody case},  arXiv: math/0703039.
  \bibitem{HU} A. Hubery, \textit{Quiver representations respecting a
      quiver automorphism: a generalisation of a theorem of Kac},  J. London Math. Soc. (2)  69  (2004),  no. 1, 79--96.
  \bibitem{Kac} V. G. Kac, \textit{Infinite dimensional Lie algebras}, Cambridge University Press, 1994
  \bibitem{L1} G. Lusztig,  \textit{Semicanonical bases
        arising from enveloping algebras},  Adv. Math.  151  (2000),  no. 2, 129--139. 
  \bibitem{L2} G. Lusztig, \textit{Quivers, perverse sheaves, and
      quantized enveloping algebras},  J. Amer. Math. Soc.  4  (1991),  no. 2, 365--421. 
  \bibitem{L3} G. Lusztig, \textit{Introduction to quantum groups},
    Progress in Mathematics 110, Birk\"auser, Boston, 1993.
  \bibitem{RR} I. Reiten, C. Riedtmann, \textit{Skew group algebras in the representation theory of artin algebras}, J. Algebra  92  (1985),  no. 1, 224--282.
\end{thebibliography}
\end{document}